\documentclass[12pt]{article}
\usepackage{graphicx,psfrag,epsfig, color,float}
\usepackage{amssymb,amsmath,amscd,amsthm}
\usepackage{graphicx,psfrag,epsfig}

\usepackage{graphicx}
\usepackage[active]{srcltx}

\newtheorem{theorem}{Theorem}[section]

\newtheorem{lemma}[theorem]{Lemma}
\newtheorem{proposition}[theorem]{Proposition}

\date{}

\begin{document}

\date{}
\title{On diffusions in media with pockets of large diffusivity}
\author{Mark Freidlin\footnote{Dept of Mathematics, University of Maryland,
College Park, MD 20742, mif@math.umd.edu}, Leonid
Koralov\footnote{Dept of Mathematics, University of Maryland,
College Park, MD 20742, koralov@math.umd.edu}, Alexander Wentzell\footnote{Dept of Mathematics, Tulane University, New Orleans, LA  70118, wentzell@math.tulane.edu}
} \maketitle

\begin{abstract}
We consider diffusion processes in media with pockets of large diffusivity. The asymptotic behavior of such processes is described when the diffusion coefficients in the pockets tend to infinity. The limiting process is identified as a diffusion on the space where each of the pockets is treated as a single point, and certain conditions on the behavior of the process on the boundary of the pockets are imposed. Calculation of various probabilities and expectations related to the limiting process leads to new initial-boundary (and boundary) problems for the corresponding parabolic (and elliptic) PDEs. 
\end{abstract}



\section{Introduction.} \label{intro}

We will study a family of diffusion processes $X^{x,\varepsilon}_t$, whose generators are of the form $M + \varepsilon^{-1} L$, where
the coefficients of $L$ are  equal to zero outside of a union of disjoint domains. Here $x$ stands for the initial position of the process and $\varepsilon$ is 
a parameter that tends to zero.  Let us assume that the state space of the processes is the $d$-dimensional torus $\mathbb{T}^d$ (general manifolds, whether compact or not, can also be considered, but we'll stick with the torus for the sake of simplicity of notation).

Let $D_1,...,D_n \subset \mathbb{T}^d $ be domains with $C^3$-smooth boundaries such that the closures $\overline{D}_1,...,\overline{D}_n$ are disjoint. Let $ U = \mathbb{T}^d \setminus \bigcup_{k=1}^n \overline{D}_k$. We assume that $U \neq  \emptyset$. Define the following operators
\begin{equation} \label{operators}
M u(x) = \frac{1}{2} \Delta u (x) ,~~~~~L u(x) = \frac{1}{2} \nabla (a(x) \cdot \nabla u (x)),
\end{equation}
where  $a \in C(\mathbb{T}^d) \bigcap C^3 (\bigcup_{k=1}^n \overline{D}_k)$ is such that $a(x) = 0$ for $x \in U$, $a(x) > 0$ for $x \in D_k$, $k =1,...,n$. 
Moreover, we'll assume that $c_1 ({\rm dist} (x, \partial D_k))^2 \leq a(x) \leq c_2 ({\rm dist} (x, \partial D_k))^2$ for some positive constants $c_1, c_2$, each $k$, and all $x \in D_k$. 

Our goal is to show that $X^{x,\varepsilon}_t$ converge, on an appropriate state space, to a limiting family of processes $Y^{x}_t$ when $\varepsilon \downarrow 0$.  First, let us give an intuitive description of the behavior of $X^{x,\varepsilon}_t$ for small $\varepsilon$. For a small $\delta > 0$, let 
\[
D^{+\delta}_k = \{x \in \mathbb{T}^d :  {\rm dist} (x, D_k) < \delta \}, ~~~ D^{-\delta}_k = \{x \in D_k : {\rm dist} (x, \partial D_k) > \delta \}.
\]
Assume that $x \in U$, in which case $X^{x,\varepsilon}_t$ moves as a Brownian motion until it reaches $\partial D_k$ for some $k$. Once the process reaches $\partial D_k$, it will reach $D^{- \delta}_k$ very soon, and will then move very fast in the interior of $D_k$ due to the large parameter $\varepsilon^{-1}$  at $L$. Next, we need to understand
how the process exits $D^{+\delta}_k$ (if $\varepsilon$ is small, the process will reach $\partial D_k$ and then go back to $D^{- \delta}_k$ many times prior to exiting $D^{+\delta}_k$). We will argue that the distribution of the exit point from $D^{+\delta}_k$ is nearly uniform with respect to the $(d-1)$-dimensional  volume measure (i.e., the measure
corresponding to the volume form  on 
$\partial D^{+\delta}_k$)  if 
$\delta$ and $\varepsilon$ are small (a small $\delta$ is fixed first, and then $\varepsilon$ is taken to zero). Thus, if we disregard the time spent inside $D_k$, the
process gets almost immediately re-distributed along $\partial D_k$ (or rather $\partial D^{+\delta}_k$) upon reaching $\partial D_k$.

Processes with somewhat similar behavior were considered in \cite{FKW1}. Instead of a self-adjoint operator multiplied by a large parameter, which amounts to fast motion inside $D_k$, the re-distribution along the boundary in \cite{FKW1} was due to a trapping mechanism: the motion inside $D_k$ was assumed to be nearly deterministic with a large vector field pointing
inside the domain. The exit times were exponentially long as $\varepsilon \downarrow 0$, and the exit was due to large deviations, while now the exit times will tend to zero. 

Due to fast mixing inside $D_k$, it is impossible to distinguish between different points of $\overline{D}_k$ without time re-scaling  when studying the limiting behavior of 
$X^{x,\varepsilon}_t$.
Let $U'$ be the metric space obtained from $\overline{U}$ by identifying all points of $\partial D_k$, turning every $\partial D_k$, $k =1,...,n$, into one point $d_k$.
The family of limiting processes $Y^x_t$, $x \in U'$,  will be defined in terms of its generator. Since we expect $Y^x_t$ to coincide with a Wiener process inside $U$, the generator coincides with $\frac{1}{2} \Delta$ on a certain class of functions. The domain of the generator of the limiting process, 
however, should be restricted by certain boundary conditions to account for non-trivial behavior of $Y^x_t$ on the boundary of $U$ and for the delay at the points $d_k$. 

The generator of the limiting process will be carefully defined in the Section~\ref{dttp}, where we also formulate the main result on convergence in the case when (\ref{operators}) holds. The  process $Y^x_t$ defined by this generator spends a positive proportion of time in $d_k$, akin to a sticky one-dimensional Brownian motion. 

The problem studied in this paper can be considered as one concerning the long time influence of a small non-degenerate perturbation ($\Delta/2$) of a degenerate diffusion (operator $L$). A related problem was studied in \cite{FH}. There, the diffusion matrix of operator $L$ is assumed to be smooth and have rank $d-1$ 
outside of the domains $D_1,...,D_n$, and full rank inside the domains. The state space for the limiting process is then a graph, whose vertices correspond to the
domains $D_k$. In the current paper, the diffusion matrix completely degenerates outside $\bigcup_{k=1}^n \overline{D}_k$, and the state space $U'$ is 
$d$-dimensional.

\section{Description of the limiting process. Formulation of the main result.} \label{dttp}

Let $X^{x,\varepsilon}_t$ be the process on $\mathbb{T}^d$ starting at $x$, with the generator $M + \varepsilon^{-1} L$. We assume that $M$ and $L$ are given 
by (\ref{operators}). Let $\varphi: \mathbb{T}^d \rightarrow U'$ be the mapping defined by $\varphi(x) = x$ for $x \in U$ and $\varphi(x) = d_k$ for $x \in \overline{D}_k$. Let $Y^{x, \varepsilon}_t = \varphi(X^{x,\varepsilon}_t)$.  Note that the superscript $x$ in the notation for the process $Y^{x, \varepsilon}_t$ is an element
of $\mathbb{T}^d$, while the process itself is $U'$-valued. 

In this section, we define the family of processes $Y^x_t$, $x \in U'$, which later will be proved to be the limiting processes for $Y^{x, \varepsilon}_t$, $x \in \mathbb{T}^d$,  as $\varepsilon \downarrow 0$. We'll use the Hille-Yosida theorem stated here in the form that is convenient for considering closures of linear operators (see \cite{We}).
\begin{theorem}
Let $K$ be a compact space, $C(K)$ be the space of continuous functions on it. The space $C(K)$ is endowed with the supremum norm.  Suppose that a linear operator $A$ on $C(K)$ has the following properties:

(a) The domain $\mathcal{D}(A)$ is dense in $C(K)$;

(b) The constant function  $ \mathbf{1}$ belongs to $\mathcal{D}(A)$ and $A \mathbf{1} = 0$;

(c) The maximum principle: If $S$ is the set of points where a function $f \in \mathcal{D}(A)$ reaches its maximum, then $A f (x) \leq 0$ for at
least one point $x \in S$.

(d) For a dense set $\Psi \subseteq C(K)$, for every $\psi \in \Psi$, and every $\lambda > 0$, there exists a solution $f \in \mathcal{D}(A)$ of the equation $\lambda f - A f = \psi$.

Then the operator $A$ is closable and its closure $ \overline{A}$ is the infinitesimal generator of a unique semi-group of positivity-preserving operators $T_t$, $t \geq 0$, on $C(K)$ with $T_t \mathbf{1} = \mathbf{1}$, $||T_t|| = 1$.
\end{theorem}

Suppose that we are given non-zero finite measures $\nu_1,...,\nu_n$ concentrated on $\partial D_1$,...,$\partial D_n$, respectively.
The Hille-Yosida theorem will be applied to the space $K = U'$, where $U'$ is the compact metric space defined above. Let us define the linear operator $A$ in $C(U')$. First we define its domain. It consists of all functions $f \in C(U')$ that satisfy the following conditions:

(1) $f$  is twice continuously differentiable in $U$;

(2) The limits of all the first and second order derivatives of $f$ exist at all the points of the boundary $\partial U = \bigcup_{k=1}^n \partial {D}_k$;

(3) There are constants $g_1,...,g_n$ such that
\[
\lim_{ x \in U, {\rm dist}(x, \partial D_k) \downarrow 0} \Delta f(x) = g_k,~~~k =1,...,n;
\]

(4)  For each $k =1,...,n$,
\begin{equation} \label{intco}
\int_{\partial D_k} \langle \nabla f(x) , n(x) \rangle \nu_k(d x) + \frac{1}{2} g_k {\rm Vol}(D_k) = 0,
\end{equation}
where  $n(x)$ is the unit exterior  normal at $x \in \partial D_k$ (external with respect to $U$).

For $f \in \mathcal{D}(A)$ and $x \in U'$, we define
\[
A f =  \begin{cases} \frac{1}{2} \Delta f(x), & {\rm if }~ x \in U, \\ \frac{1}{2} g_k, & {\rm if }~ x = d_k,~~~k=1,...,n. \end{cases}
\]
Let us check that the conditions of the Hille-Yosida theorem are satisfied.

(a) Consider the set $G$ of functions $g$ on $U'$ that are infinitely differentiable on $U$ and have the following property: for each $k =1,...,n$ there is
a set $V_k$ open in $U'$ such that $d_k \in V_k$ and $g$ is constant on $V_k$. It is clear that $G \subset \mathcal{D}(A)$ and $G$ is dense in $C(U')$.

(b) Clearly $\mathbf{1} \in \mathcal{D}(A)$ and $A \mathbf{1} = 0$.

(c) If $f$ has a maximum at $x \in U$, it is clear that $ \Delta f (x) \leq 0$. Now suppose that $f$ has a maximum at $d_k$. We can view $f$ as an element of  $C^2(\overline{U})$ that is constant on each component of the boundary, in particular on $\partial D_k$. Let $g_k$ be the value of $\Delta f$ on $\partial D_k$, and suppose that $g_k > 0$. Then, by (\ref{intco}),  $\langle \nabla f(x) , n(x) \rangle$ must be negative for some $x \in \partial D_k$, which contradicts the fact that $f$ reaches its maximum on $\partial D_k$. Therefore, $g_k \leq 0$, as required. 

(d) Let $\Psi$ be the set of functions $\psi \in C(U')$ that have limits of all the first order derivatives as $y \rightarrow x, y \in U$, at all points
$x \in \partial U$. It is clear that $\Psi$ is dense in $C(U')$. Let $\widetilde{f} \in C^2 (\overline{U})$ be the solution of the equation
$\lambda \widetilde{f} - \frac{1}{2} \Delta \widetilde{f} = \psi$ in $U$, $\widetilde{f} = 0$ on $\partial U$. Let $h_k \in C^2(\overline{U})$ be the solution of the equation
\[
\lambda {h}_k(x) - \frac{1}{2} \Delta {h}_k(x) = 0,~~~~x \in U,
\]
\[
h_k(x) = 1,~~x \in \partial D_k;~~~~~h_k(x) = 0,~~x \in \partial U \setminus \partial D_k.
\]
Let us look for the solution $f \in \mathcal{D}(A)$ of $\lambda f -  A f = \psi$ in the form $f = \widetilde{f} + \sum_{k =1}^n c_k h_k$. We get $n$ linear equations for $c_1,...,c_n$. The solution is unique because of the maximum principle. Therefore, the determinant of the system is non-zero, and the solution exists for all the right hand sides.

Let $ \overline{A}$ be the closure of $A$. Let $T_t$, $t \geq 0$, be the corresponding semi-group on $C(U')$, whose existence is guaranteed by the
Hille-Yosida theorem. By the Riesz-Markov-Kakutani representation theorem, for $x \in U'$ there is a measure $P(t,x,dy)$ on $(U', \mathcal{B}(U'))$ such that
\[
(T_t f)(x) = \int_{U'} f(y)P(t,x,dy),~~~f \in C(U').
\]
It is a probability measure since $T_t \mathbf{1} = \mathbf{1}$. Moreover, it can be easily verified that $P(t,x,B)$ is a Markov transition function.  Let $Y^x_t$, $x \in U'$, be the corresponding Markov family. In order to show that a modification with continuous trajectories exists, it is enough to check that $\lim_{t \downarrow 0} P(t,x,B)/t = 0$ for each closed set $B$ that doesn't contain $x$  (see \cite{Dyn} or Theorem I.5 of \cite{Mandl}).  Let  $f \in \mathcal{D}(A)$  be a non-negative function that is equal to one on $B$ and is equal to zero in some neighborhood of $x$. Then
\[
\lim_{t \downarrow 0} \frac{P(t,x,B)}{t} \leq \lim_{t \downarrow 0} \frac{(T_t f)(x) - f(x)}{t} = Af(x) = 0,
\]
as required. Thus $Y^x_t$ can be assumed to have continuous trajectories.

Now we are prepared to formulate the main result on the convergence of the processes $X^{x, \varepsilon}_t$. Let $\nu_1,...,\nu_n$ be the $(d-1)$-dimensional volume measures (i.e., the measures corresponding to the volume forms) on $\partial D_1,...,\partial D_n$, respectively. Let  $Y^{\varphi(x)}_t$ be the corresponding process, constructed above, with values in $U'$.

\begin{theorem} \label{mt1a} Let $X^{x,\varepsilon}_t$ be the process on $\mathbb{T}^d$ starting at $x$, with the generator $M + \varepsilon^{-1} L$.
For each $x \in \mathbb{T}^d$, the measures on $C([0,\infty), U')$ induced by the processes $Y^{x, \varepsilon}_t = \varphi(X^{x,\varepsilon}_t)$ converge weakly, 
as $\varepsilon \downarrow 0$, to the measure induced by $Y^{\varphi(x)}_t$.
\end{theorem}

\section{The theorem on the convergence of the processes} \label{seconv}
In this section, we prove Theorem~\ref{mt1a}. We'll use the convention that a function $f$ defined on $U'$ can also be viewed as a function on $\mathbb{T}^d$. The key ingredient in the proof of Theorem~\ref{mt1a} is the following proposition.
\begin{proposition} \label{mainpro}
Suppose that $f \in \mathcal{D}(A)$. Then
\begin{equation} \label{maeq}
\lim_{\varepsilon \downarrow 0} \mathrm{E} \left(f(Y^{x, \varepsilon}_t) - f(x) -  \frac{1}{2} \int_0^t \Delta f (Y^{x, \varepsilon}_u) du \right) = 0
\end{equation}
for each $t \geq 0$, uniformly in  $x \in \mathbb{T}^d$.
\end{proposition}
\proof 
For the sake of notational simplicity, we'll assume that there is just one domain where the coefficients of $L$ are non-zero. The proof in the case of multiple domains is similar. We'll denote the domain by $D$ and will drop the subscript $k$ from the notation everywhere. For example, (\ref{intco}) now takes the form
\begin{equation} \label{ortt}
\int_{\partial D} {\langle \nabla f (x), n(x) \rangle}   \nu(dx)  + \frac{1}{2} g {\rm Vol}(D) = 0.
\end{equation}

The proof will rely on several lemmas concerning the behavior of the process $X^{x,\varepsilon}_t$ in the vicinity of $D$. We'll state these lemmas when they are needed,
while the proofs will be provided after we complete the proof of the proposition.

Fix an arbitrary $\eta > 0$. Our goal is to show that
\[
\Big|\mathrm{E} \left(f(Y^{x, \varepsilon}_t) - f(x) -  \frac{1}{2} \int_0^t \Delta f (Y^{x, \varepsilon}_u) du \right)  \Big| \leq \eta 
\]
for all sufficiently small $\varepsilon$, uniformly in  $x \in \mathbb{T}^d$.  For $B \subset \mathbb{T}^d$, let 
\[
\tau^{x, \varepsilon}(B) =  \inf \{t\geq 0: X^{x,\varepsilon}_t \in B \}.
\]
We'll use the following lemma, which will be proved in the end of this section.
\begin{lemma} \label{jkjl}
The following limit
\[
\lim_{\varepsilon \downarrow 0} \mathrm{E} \tau^{x, \varepsilon}(\partial D) = 0
\] 
is uniform in $x \in D$. 
\end{lemma}

For $\delta > 0$, we define two sequences of stopping times:
\[
\sigma^{x, \varepsilon}_0 = 0,
\]
\[
\tau^{x, \varepsilon}_n = \inf(t \geq \sigma^{x, \varepsilon}_{n-1}: X^{x, \varepsilon}_t \in \partial D),~~ n \geq 1,
\]
\[  
\sigma^{x, \varepsilon}_n = \inf(t \geq \tau^{x, \varepsilon}_n: X^{x, \varepsilon}_t \in \partial D^{+\delta}),~~n \geq 1. 
\]
Then, extending $f$ from $U'$ to $\mathbb{T}^d$ (as a constant on $\overline{D}$ equal to the value of the original function at $d$) and putting $\Delta f \equiv g $ on $D$, we obtain
\[
\mathrm{E} \left(f(Y^{x, \varepsilon}_t) - f(x) -  \frac{1}{2}  \int_0^t \Delta f (Y^{x, \varepsilon}_u) du \right) =
\]
\[
\mathrm{E} \sum_{n = 1}^\infty  \left(f(X^{x, \varepsilon}_{\tau^{x, \varepsilon}_n \wedge t} ) - f(X^{x, \varepsilon}_{\sigma^{x, \varepsilon}_{n-1} \wedge t}) - \frac{1}{2}  \int_{\sigma^{x, \varepsilon}_{n-1} \wedge t}^{\tau^{x, \varepsilon}_n \wedge t} \Delta f (X^{x, \varepsilon}_u) du \right) +
\]
\[
\mathrm{E} \sum_{n = 1}^\infty  \left(f(X^{x, \varepsilon}_{\sigma^{x, \varepsilon}_n \wedge t} ) - f(X^{x, \varepsilon}_{\tau^{x, \varepsilon}_{n} \wedge t}) - \frac{1}{2}  \int_{\tau^{x, \varepsilon}_{n} \wedge t}^{\sigma^{x, \varepsilon}_n \wedge t} \Delta f (X^{x, \varepsilon}_u) du \right).
\]
The first expectation on the right hand side tends to zero  uniformly in $x \in \mathbb{T}^d$. Indeed, if $x \in  \overline{U}$, then all the terms are equal to zero since $X^{x, \varepsilon}_t$ is a Wiener process on $\overline{U}$. If $x \in D$, then   all the terms but the first one are equal to zero. The first term is equal to 
$ {g} \mathrm{E} \tau^{x, \varepsilon}(\partial D)/2$,
which tends to zero by Lemma~\ref{jkjl} uniformly in $x \in D$.

Our next goal is to select a sufficiently small $\delta$ so the second expectation can be made smaller than $2\eta/3$ for all sufficiently small $\varepsilon$.
\\
\\ 
{\bf Remark on notation.} From this point on, we'll use the notation ${\alpha}^{\varepsilon}(\delta)$ for a generic quantity that satisfies  
\[
\lim_{\delta \downarrow 0} \limsup_{\varepsilon \downarrow 0}  |{\alpha}^{\varepsilon}(\delta)| = 0.
\] 
The notation may stand for different quantities from line to line. 
\\

 We will use the
following lemma, also to be proved in the end of this section. 
\begin{lemma} \label{exde}
The following asymptotic  relations hold  true: 
\[
\sup_{x \in D} | 
\mathrm{E} \tau^{x, \varepsilon}(\partial D^{+\delta}) - \delta \frac{  {\rm Vol} (D) }{  \nu(\partial D)} | =  \delta {\alpha}^{\varepsilon}(\delta).
\] 
\[
\sup_{x \in D^{+\delta}} | 
\mathrm{E} \tau^{x, \varepsilon}(\partial D^{+\delta}) | =  {\alpha}^{\varepsilon}(\delta).
\]
\end{lemma}

Let $N = N(x, \delta, \varepsilon) = \max\{n: \tau^{x, \varepsilon}_{n} < t \}$. Observe that
\[
\sup_{x \in \mathbb{T}^d} \Big| \mathrm{E} \sum_{n = 1}^\infty  \left(f(X^{x, \varepsilon}_{\sigma^{x, \varepsilon}_n \wedge t} ) - f(X^{x, \varepsilon}_{\tau^{x, \varepsilon}_{n} \wedge t}) - \frac{1}{2}  \int_{\tau^{x, \varepsilon}_{n} \wedge t}^{\sigma^{x, \varepsilon}_n \wedge t} \Delta f (X^{x, \varepsilon}_u) du \right) - 
\]
\[
\mathrm{E} 
\sum_{n = 1}^N\left(f(X^{x, \varepsilon}_{\sigma^{x, \varepsilon}_n } ) - f(X^{x, \varepsilon}_{\tau^{x, \varepsilon}_{n} }) - \frac{1}{2}  
\int_{\tau^{x, \varepsilon}_{n} }^{\sigma^{x, \varepsilon}_n } \Delta f (X^{x, \varepsilon}_u) du \right)  \Big| \leq
\]
\[
\sup_{x \in D^{+\delta}} \Big| \mathrm{E} \left(f(X^{x, \varepsilon}_{\tau^{x, \varepsilon}(\partial D^{+\delta}) } ) - f(x) - \frac{1}{2}  
\int_{0 }^{\tau^{x, \varepsilon}(\partial D^{+\delta}) } \Delta f (X^{x, \varepsilon}_u) du \right) \Big| \leq
\]
\[
\delta \sup_{x \in \overline{U}} |\nabla f(x)  | + \frac{1}{2}{\sup_{x \in \overline{U}} |\Delta f(x)  |} \sup_{x \in D^{+\delta}} \mathrm{E} \tau^{x, \varepsilon}(\partial D^{+\delta}).
\]
By Lemma~\ref{exde}, there is $\delta_0 > 0$ and for each $0 \leq \delta < \delta_0$ there is $\varepsilon_1(\delta)$ such that
\[
\delta \sup_{x \in \overline{U}} |\nabla f(x)  | + \frac{1}{2}{\sup_{x \in \overline{U}} |\Delta f(x)  |} \sup_{x \in D^{+\delta}} \mathrm{E} \tau^{x, \varepsilon}(\partial D^{+\delta}) < \frac{\eta}{3}
\]
provided that  $\varepsilon < \varepsilon_1(\delta)$. Therefore it remains to find a sufficiently small $\delta$ so that
\[
\Big| \mathrm{E} 
\sum_{n = 1}^N\left(f(X^{x, \varepsilon}_{\sigma^{x, \varepsilon}_n } ) - f(X^{x, \varepsilon}_{\tau^{x, \varepsilon}_{n} }) - \frac{1}{2}  
\int_{\tau^{x, \varepsilon}_{n} }^{\sigma^{x, \varepsilon}_n } \Delta f (X^{x, \varepsilon}_u) du \right)  \Big| \leq \frac{\eta}{3}
\]
for all sufficiently small $\varepsilon$. We start by noting that
\[
\sup_{x \in \mathbb{T}^d} \Big| \mathrm{E} 
\sum_{n = 1}^N\left(f(X^{x, \varepsilon}_{\sigma^{x, \varepsilon}_n } ) - f(X^{x, \varepsilon}_{\tau^{x, \varepsilon}_{n} }) - \frac{1}{2}  
\int_{\tau^{x, \varepsilon}_{n} }^{\sigma^{x, \varepsilon}_n } \Delta f (X^{x, \varepsilon}_u) du \right)  \Big| \leq 
\]
\begin{equation} \label{toyy1}
(\sup_{x \in \mathbb{T}^d} \mathrm{E} N ) \sup_{x \in \partial D} \Big| \mathrm{E} \left(f(X^{x, \varepsilon}_{\tau^{x, \varepsilon}(\partial D^{+\delta}) } ) - f(x) - \frac{1}{2}  
\int_{0 }^{\tau^{x, \varepsilon}(\partial D^{+\delta}) } \Delta f (X^{x, \varepsilon}_u) du \right) \Big|.
\end{equation}
Recall, from the definition of $N$, that the process $X^{x, \varepsilon}_t$ makes no fewer than $N -1$ excursions from $\partial D^{+\delta}$ to $\partial D$ during the time interval $[0,t]$. Since the process coincides with the Brownian motion on $\overline{U}$, there is a constant $c = c(t)$ such that
\begin{equation} \label{ey1}
\sup_{x \in \mathbb{T}^d} \mathrm{E} N \leq c \delta^{-1}
\end{equation}
for all sufficiently small $\delta$ and for all $\varepsilon$. 

For $x \in \partial D^{+\delta}$, we define $\theta (x) \in \partial D$ as the point such that ${\rm dist}(x, \theta(x)) = {\rm dist}(x, \partial D)$.
Since $\partial D$ is smooth, $\theta(x)$ is defined uniquely for each $x$, provided that $\delta$ is sufficiently small. 

Let $\mu^{x, \varepsilon}_\delta$ be the distribution of $X^{x, \varepsilon}_{\tau^{x, \varepsilon}(\partial D^{+\delta})}$  (it's a probability measure on
$\partial D^{+\delta}$). We claim that it is close to the normalized $(d-1)$-dimensional volume measure $\overline{\nu} = \nu/\nu(\partial D)$ on $\partial D$. More precisely, we have the following lemma, which will be proved in the end of this section.
\begin{lemma} \label{le3rr}
Let $\varphi$ be a continous function on $\partial D$.  Then   
\[
\sup_{x \in \partial D} \Big|   \int_{\partial D^{+\delta}} \varphi(\theta(y))  \mu^{x, \varepsilon}_\delta(dy)  - \int_{\partial D} \varphi (y) 
\overline{\nu} (dy)  \Big| 
= \alpha^\varepsilon(\delta). 
\]
\end{lemma}

For $y \in  \partial D^{+\delta}$, we have, by the Taylor formula,
\[
f(y) = f(\theta(y)) - \langle \nabla f(\theta(y)) , n(\theta(y)) \rangle \delta + r(y, \delta) \delta^2, 
\]
where $n(x)$ is the exterior normal (with respect to $U$) at $x$ and $r$ is a bounded function. Applying this 
to $y = X^{x, \varepsilon}_{\tau^{x, \varepsilon}(\partial D^{+\delta}) }$ and using Lemma~\ref{le3rr} with $\varphi = \langle 
\nabla f , n \rangle$, we obtain 
\begin{equation} \label{ey2}
\sup_{x \in \partial D} \Big| \mathrm{E} \left(f(X^{x, \varepsilon}_{\tau^{x, \varepsilon}(\partial D^{+\delta}) } ) - f(x) + \delta \int_{\partial D}  \langle 
\nabla f , n \rangle (y) \overline{\nu}(dy)  \right) \Big| = \delta \alpha^\varepsilon(\delta). 
\end{equation}
By Lemma~\ref{exde}, there is $c > 0$ such that 
\[
\sup_{x \in \partial D} \Big| \mathrm{E} \left(\frac{1}{2}  
\int_{0 }^{\tau^{x, \varepsilon}(\partial D^{+\delta}) } \Delta f (X^{x, \varepsilon}_u) du  - \frac{ g {\rm Vol}(D)  \delta }{2 \nu(\partial D)} \right) \Big| \leq 
\]
\begin{equation} \label{ey3}
 c \delta (\alpha^\varepsilon(\delta) + \sup_{x \in  D^{+\delta}}|\Delta f(x) - g|) = \delta \alpha^\varepsilon(\delta).
\end{equation}
 Combining  (\ref{ey1}),
(\ref{ey2}), and (\ref{ey3}), we obtain an estimate on the right hand side of (\ref{toyy1}): 
\[
(\sup_{x \in \mathbb{T}^d} \mathrm{E} N ) \sup_{x \in \partial D} \Big| \mathrm{E} \left(f(X^{x, \varepsilon}_{\tau^{x, \varepsilon}(\partial D^{+\delta}) } ) - f(x) - \frac{1}{2}  
\int_{0 }^{\tau^{x, \varepsilon}(\partial D^{+\delta}) } \Delta f (X^{x, \varepsilon}_u) du \right) \Big| \leq 
\]
\[
 c \delta^{-1} \left( \delta \big| \int_{\partial D}  \langle 
\nabla f , n \rangle (y) \overline{\nu} (dy) + \frac{g {\rm Vol}(D)}{2 \nu(\partial D)}  \Big|  + \delta \alpha^\varepsilon(\delta) \right).
\]
The first term in the brackets on the right hand side is equal to zero by (\ref{ortt}). Therefore, the whole expression can be made smaller than $\eta/3$ by 
selecting a sufficiently small $\delta$. This completes the proof of the proposition.
\qed
\\
\\ 
{\it Proof of Theorem~\ref{mt1a}.} Recall that $\Psi$ is the set of functions $\psi \in C(U')$ that have limits of all the first order 
derivatives as $ y \rightarrow x, y \in U$, at all points $x \in \partial U$. This is a measure-defining class of functions 
on $U'$, i.e., if $\mu_1$ and $\mu_2$ satisfy  $\int_{U'} \psi (x) \mu_1(dx) = \int_{U'}\psi(x) \mu_2(dx)$ for every $\psi \in \Psi$, 
then $\mu_1 = \mu_2$.  As shown in Section~\ref{dttp}, for every $\psi \in \Psi$ and every $\lambda > 0$, there 
is $f \in \mathcal{D}(A)$ that satisfies $\lambda f - Af = \psi$.  We have demonstrated that (\ref{maeq}) holds 
for $f \in \mathcal{D}(A)$. By Lemma 3.1. in Chapter 8 of \cite{FW}, this is sufficient to guarantee the convergence if, 
in addition, the family  $\{Y^{x, \varepsilon}_t\}, \varepsilon > 0, x \in {U'}$ is tight. The tightness, however, 
is clear since the processes coincide with a Wiener process inside $U$, while all the points of $\partial U$  and $D$  are identified.
\qed 
\\

It remains to prove Lemmas~\ref{jkjl}, \ref{exde}, and \ref{le3rr}. 
The proofs  of these lemmas will make use of the strong Doeblin condition, which we formulate next. 
Let $Y^x_t$, $x \in S$, $t \in \mathbb{R}_+$ or $t \in \mathbb{Z}_+$, be a Markov family on 
$S \in \mathcal{B}(\mathbb{R}^d)$. It is said to satisfy the strong Doeblin condition if there are a probability measure  $Q$  on $\mathcal{B}(S)$,
a constant $c > 0$,  and a time $t_0$ such that 
\begin{equation} \label{doco}
\mathrm{P}(Y^x_{t_0} \in A) \geq c  Q(A)
\end{equation} 
 for all $x \in S$, $A \in  \mathcal{B}(S)$. 
\begin{lemma} \label{sdc}  {\rm (\cite{Doob}, Ch 5, 6)} 
If the strong Doeblin condition is satisfied, then there is a unique invariant measure $\mu$ for the family $Y^x_t$. Moreover, there is
$\lambda > 0$ (that can be defined in terms of $t_0$ and $c$) such that
\[
|\mathrm{P}(Y^x_{t} \in A) - \mu(A) | \leq \lambda^{-1} e^{-\lambda t},~~{\rm for}~~t \geq 0,~x \in S,~A \in \mathcal{B}(S).
\]
\end{lemma} 
The next lemma is a form of the Strong Law of Large Numbers. Its proof is similar to the proof of the  Law of Large Numbers for discrete-time Markov chains that
can be found in \cite{Doob}, Chapter 5.
\begin{lemma} \label{llnu}
For each $A \in \mathcal{B}(S)$ with $\mu(A) > 0$ and each $\varepsilon > 0$, there is $T > 0$ such that
\[
\mathrm{P}(\sup_{t \geq T} |\frac{\int_0^t \chi_A(Y^x_{s}) ds}{t \mu(A)} -1 | > \varepsilon) < \varepsilon,~~~~x \in S. 
\]
\end{lemma}
For fixed $c$, $t_0$, $A$, and $\mu$, the same  $T >0$ can be chosen for all the Markov families that satisfy (\ref{doco}) and
have $\mu$ as the invariant measure.  
\\

\noindent
{\it Proof of Lemma~\ref{jkjl}.} Let 
\[
G^\delta =  \overline{D} \setminus D^{-\delta} = \{x \in \overline{D} : {\rm dist} (x, \partial D) \leq \delta \}.
\]
For $x \in G^\delta$, let $h(x) = {\rm dist} (x, \partial D)$. Recall that $\theta (x) \in \partial D$ is the point such that 
${\rm dist}(x, \theta(x)) = {\rm dist}(x, \partial D)$. If $\delta$ is small enough, then $x \rightarrow (h, \theta)$ is a $C^3$-diffeomorphism
between $G^\delta$ and $[0,\delta] \times \partial D$. First, let us provide an upper estimate on the function
\[
g^\delta_\varepsilon (x) = \mathrm{E} \int_0^{\tau^{x, \varepsilon}(\partial D)} \chi_{G^\delta} (X^{x, \varepsilon}_t) dt,~~~~x \in  \overline{D}, 
\]
which is the expectation of the time the process spends in $G^\delta$ prior to exiting $D$.  

From our assumptions on the function $a$, it follows that there is a positive continuous function $\psi \in C(\partial D)$
 such that 
\begin{equation} \label{qq1}
a(x) = \psi(\theta(x)) h^2(x)+ K_1  (x)  h^3(x), 
\end{equation}
\[
\nabla a (x)  = 2 \psi(\theta(x)) h(x) \nabla h(x)  + K_2  (x)  h^2(x),
\]
where $K_1$ and $K_2$ are bounded on $G^\delta$  ($K_2$ is vector-valued). 
Define 
\[
r =  \inf_{x \in \partial D} \psi(x),~~~R = \sup_{x \in \partial D} \psi(x).
\]
Let
\[
w^\delta_\varepsilon(h) =  \frac{2R}{r}\int_0^h \frac{\delta -t}{\varepsilon + R t^2} dt ,~~~0 \leq h \leq \delta,
\]
and 
\[
u^\delta_\varepsilon(x) = \begin{cases} 2\varepsilon w^\delta_\varepsilon(h(x)), & \mbox{if  } x \in G^\delta \\ 2 \varepsilon w^\delta_\varepsilon(\delta), & \mbox{if  } x \in D \setminus G^\delta. \end{cases}
\]

Observe that there is $C > 0$ such that $
u^\delta_\varepsilon (x) \leq C \sqrt{\varepsilon}$ for $x \in \overline{D}$.
 We will show that for all sufficiently small $\delta$ we have
\begin{equation} \label{enm1}
(M + \varepsilon^{-1} L) u^\delta_\varepsilon (x) \leq -1,~~~x \in {\rm Int}(G^\delta). 
\end{equation}
Observe that  $(w^\delta_\varepsilon)'(\delta) = 0$. Therefore, the first order partial derivatives of $u^\delta_\varepsilon$
are Lipschitz continuous on $D$ (including $\partial D^{-\delta}$).  Also note that the funcion $u^\delta_\varepsilon$ is $C^2$-smooth on $\overline{D} \setminus \partial D^{-\delta}$,
while the second derivatives may have a 
jump at the surface $\partial D^{-\delta}$. For $x \in \overline{D}$, we can apply the Ito formula to $u^\delta_\varepsilon(X^{x,\varepsilon}_t)$
on the time interval $[0,\tau^{x, \varepsilon}(\partial D)]$ (this
is justified by approximating the function $u^\delta_\varepsilon$, for fixed $\delta$ and $\varepsilon$, by a sequence of twice
continuously differentiable functions whose second derivatives are uniformly bounded). Thus
\[
\mathrm{E} u^\delta_\varepsilon (X^{x,\varepsilon}_{ \tau^{x, \varepsilon}(\partial D)}) =  
u^\delta_\varepsilon (x) + \mathrm{E} \int_0^{ \tau^{x, \varepsilon}(\partial D)}  (M + \varepsilon^{-1} L) u^\delta_\varepsilon (X^{x,\varepsilon}_{t}) dt.
\]
Note that the left hand side here is equal to zero since $u^\delta_\varepsilon(x) = w^\delta_\varepsilon(0) = 0$ for $x \in \partial D$. Therefore,
\[
u^\delta_\varepsilon (x) = - \mathrm{E} \int_0^{ \tau^{x, \varepsilon}(\partial D)}  (M + \varepsilon^{-1} L) u^\delta_\varepsilon (X^{x,\varepsilon}_{t}) dt \geq
 \mathrm{E} \int_0^{\tau^{x, \varepsilon}(\partial D)} \chi_{G^\delta} (X^{x, \varepsilon}_t) dt = g^\delta_\varepsilon (x),~~~x \in \overline{D}. 
\]

 In order to verify (\ref{enm1}), we write
\[
(M + \varepsilon^{-1} L) u^\delta_\varepsilon (x) =  (\varepsilon +  a(x)) (w^\delta_\varepsilon)'' (h(x)) +   
\]
\[
 \left(  \langle \nabla a, \nabla h \rangle + (\varepsilon + a(x)) \Delta h \right) (w^\delta_\varepsilon)' (h(x)) = A_1 + A_2 + A_3, 
\]
where
\[
A_1 =   (\varepsilon + \psi(\theta(x)) h^2(x) )  (w^\delta_\varepsilon)'' (h(x)) +  2\psi(\theta(x)) h(x) (w^\delta_\varepsilon)' (h(x)),
\]
\[
A_2 =  K_1 h^3 (x)  (w^\delta_\varepsilon)'' (h(x)) ,
\]
\[
A_3 = \left(  \langle K_2, \nabla h \rangle h^2 (x)  + (\varepsilon + a(x)) \Delta h \right) (w^\delta_\varepsilon)' (h(x)).
\]
From the definition of $w^\delta_\varepsilon$ it follows that $A_1 \leq -2$. From (\ref{qq1}) and the boundedness of $K_1$, $\|K_2\|$, 
it easily follows that $|A_2|, |A_3| \leq 1/2$ when $x \in G^\delta$, provided that $\delta$ is sufficiently small. Thus (\ref{enm1}) holds, and
therefore there is a constant $C = C(\delta)$ such that
\begin{equation} \label{jjip}
\mathrm{E} \int_0^{\tau^{x, \varepsilon}(\partial D)} \chi_{G^\delta} (X^{x, \varepsilon}_t) dt = g^\delta_\varepsilon (x) \leq u^\delta_\varepsilon (x) \leq C \sqrt{\varepsilon},~~~x \in \overline{D}.
\end{equation}

Next, with $\delta$ sufficiently small for this inequality to hold, 
consider the process $\widehat{X}^{x, \varepsilon}_t$ obtained from ${X}^{x, \varepsilon}_t$ by first slowing it down using time change with the 
factor $\varepsilon$ and then running the clock only when the process is inside  
$D \setminus {\rm Int}({G^{\delta/2}}) = \overline{D^{-\delta/2}}$. More precisely, let
\[
 \beta(t)  = \int_0^{t} \chi_{ D \setminus {\rm Int}({G^{\delta/2}}) } (X^{x, \varepsilon}_{\varepsilon u}) du,
\]
\[
\eta(t) = \inf \{s:  \beta(s)  > t \},
\]
and
\[
\widehat{X}^{x, \varepsilon}_t = {X}^{x, \varepsilon}_{\varepsilon \eta(t)} \in  D \setminus {\rm Int}({G^{\delta/2}}).
\]  
Since $a > 0$ on $D \setminus {\rm Int}({G^{\delta/2}})$,  the process  $\widehat{X}^{x, \varepsilon}_t$, $x \in D \setminus {\rm Int}({G^{\delta/2}})$, satisfies the  strong Doeblin condition
on $D \setminus {\rm Int}({G^{\delta/2}})$, uniformly in $\varepsilon > 0$ (a similar statement can be found in \cite{FR}, Lemma 3.7.1). The uniformity in $\varepsilon$ means that the constants $c$ and $t_0$ can be found such that (\ref{doco}), applied to the process $\widehat{X}^{x, \varepsilon}_t$, holds for all $\varepsilon$.   Therefore,  since the Lebesgue measure on $D \setminus {\rm Int}({G^{\delta/2}})$ is invariant for  $\widehat{X}^{x, \varepsilon}_t$,
 from Lemma~\ref{llnu} it follows that there is $T > 0$ such that for every stopping time $\tau \geq T$ with $\mathrm{E} \tau < \infty$ we have 
\begin{equation} \label{ner}
\mathrm{P} \left( \int_0^{\tau} \chi_{ D \setminus {\rm Int}({G^{\delta}}) } (\widehat{X}^{x, \varepsilon}_{t}) dt > 
 C \int_0^{\tau} \chi_{ G^\delta \setminus {\rm Int}({G^{\delta/2}}) } (\widehat{X}^{x, \varepsilon}_{t}) dt \right) \leq \frac{1}{4},
\end{equation}
where $C = {2 {\rm Vol} (D \setminus {\rm Int}({G^{\delta}}) ) }/{{\rm  Vol}  (  G^\delta \setminus {\rm Int}({G^{\delta/2}}))}$. Take $k > 0$, which will be specified later. For $x \in \overline{D}$, we have
\[
\mathrm{P} (\tau^{x, \varepsilon}(\partial D)  \geq k \sqrt{\varepsilon}) \leq \mathrm{P}(  \int_0^{\tau^{x, \varepsilon}(\partial D)} \chi_{G^\delta} (X^{x, \varepsilon}_t) dt \geq \frac{k \sqrt{\varepsilon}}{2}) + 
\]
\[
\mathrm{P}(  \int_0^{\tau^{x, \varepsilon}(\partial D)} \chi_{D \setminus {\rm Int} (G^\delta)} (X^{x, \varepsilon}_t) dt \geq \frac{k \sqrt{\varepsilon}}{2};~~
\tau^{x, \varepsilon}(\partial D)  \geq k \sqrt{\varepsilon} ).
\]
The first term on the right hand side is estimated from above by $1/4$ for large enough $k$, using 
 (\ref{jjip}) and the Chebyshev inequality. For the second term, we write
\[
\mathrm{P} \left(  \int_0^{\tau^{x, \varepsilon}(\partial D)} \chi_{D \setminus {\rm Int} (G^\delta)} (X^{x, \varepsilon}_t) dt \geq \frac{k \sqrt{\varepsilon}}{2};~~
\tau^{x, \varepsilon}(\partial D)  \geq k \sqrt{\varepsilon} \right) \leq 
\]
\[
\mathrm{P} \left(  \int_0^{\tau^{x, \varepsilon}(\partial D)} \chi_{D \setminus {\rm Int} (G^\delta)} (X^{x, \varepsilon}_t) dt \geq c
 \int_0^{\tau^{x, \varepsilon}(\partial D)}
 \chi_{ G^\delta \setminus {\rm Int}({G^{\delta/2}}) } ({X}^{x, \varepsilon}_{t}) dt  ;~~
\tau^{x, \varepsilon}(\partial D)  \geq k \sqrt{\varepsilon} \right) + 
\]
\[
\mathrm{P} \left( c
 \int_0^{\tau^{x, \varepsilon}(\partial D)}
 \chi_{ G^\delta \setminus {\rm Int}({G^{\delta/2}}) } ({X}^{x, \varepsilon}_{t}) dt \geq \frac{k \sqrt{\varepsilon}}{2} \right).
\]
The first term on the right hand side is estimated from above by $1/4$ for large enough $k$ using (\ref{ner}), while the second term is estimated
by $1/4$ using  (\ref{jjip}) and the Chebyshev inequality. We conclude that
\[
\mathrm{P} (\tau^{x, \varepsilon}(\partial D)  \geq k \sqrt{\varepsilon})  \leq \frac{3}{4}.
\]
Since $x \in \overline{D}$ was arbitrary, the statement of the lemma follows from here and the Markov property of the process. 
\qed
\\
\\
{\it Proof of Lemma~\ref{exde}.} First, we claim that
\begin{equation} \label{ff1}
\lim_{\varepsilon \downarrow 0} \sup_{x \in G^\delta} \mathrm{E} \tau^{x, \varepsilon}(\partial D^{+\delta} \cup \partial D^{-\delta})  = 0,
\end{equation}
 for all sufficiently small $\delta > 0$. The proof of this statement is similar to that of the previous lemma, so we omit it here. 

Define two sequences of stopping times:
\[
\widetilde{\sigma}^{x, \varepsilon}_0 = 0,
\]
\[
\widetilde{\tau}^{x, \varepsilon}_n = \inf(t \geq \widetilde{\sigma}^{x, \varepsilon}_{n-1}:  {X}^{x, \varepsilon}_t \in \partial D),~~ n \geq 1,
\]
\[  
\widetilde{\sigma}^{x, \varepsilon}_n = \inf(t \geq \widetilde{\tau}^{x, \varepsilon}_n: {X}^{x, \varepsilon}_t  \in \partial D^{- \delta}),~~n \geq 1. 
\]
Let us show that
\begin{equation} \label{cdc}
\lim_{\varepsilon \downarrow 0} \frac{\mathrm{E} \tau^{x, \varepsilon}(\partial D^{+\delta})}{\int_{\partial D^{-\delta}} \mathrm{E} \tau^{y, \varepsilon}(\partial D^{+\delta}) \widetilde{\mu}^\varepsilon(d y)} = 1
\end{equation}
uniformly in $x \in \partial D$, where $\widetilde{\mu}^\varepsilon$ is the invariant measure for the Markov 
chain ${X}^{x, \varepsilon}_{ \widetilde{\sigma}^{x, \varepsilon}_n}$.

 Let 
\[
I(\varepsilon) = \sup_{x \in \partial D^{-\delta}} | \mathrm{E} \tau^{x, \varepsilon}(\partial D^{+\delta}) -
 \int_{\partial D^{-\delta}} \mathrm{E} \tau^{y, \varepsilon}(\partial D^{+\delta})  \widetilde{\mu}^\varepsilon(d y) |,
\]
\[
J(\varepsilon) = \int_{\partial D^{-\delta}} \mathrm{E} \tau^{y, \varepsilon}(\partial D^{+\delta}) \widetilde{\mu}^\varepsilon(dy). 
\]
Observe that for each $n$
and $x \in \partial D^{-\delta}$, 
\[
\mathrm{E} \tau^{x, \varepsilon}(\partial D^{+\delta}) -
 \int_{\partial D^{-\delta}} \mathrm{E} \tau^{y, \varepsilon}(\partial D^{+\delta}) \widetilde{\mu}^\varepsilon(dy) = 
\]
\begin{equation} \label{miqu}
\left( \mathrm{E}   \tau^{ {X}^{x, \varepsilon}_{ \widetilde{\sigma}^{x, \varepsilon}_n},  \varepsilon}(\partial D^{+\delta})  -
\int_{\partial D^{-\delta}} \mathrm{E} \tau^{y, \varepsilon}(\partial D^{+\delta}) \widetilde{\mu}^\varepsilon(dy)  \right) +
\end{equation}
\[
\left( \mathrm{E} \tau^{x, \varepsilon}(\partial D^{+\delta}) - \mathrm{E}   \tau^{ {X}^{x, \varepsilon}_{ \widetilde{\sigma}^{x, \varepsilon}_n},  
\varepsilon}(\partial D^{+\delta})  \right). 
\]
The discrete-time Markov chain  
${X}^{x, \varepsilon}_{ \widetilde{\sigma}^{x, \varepsilon}_n} \in \partial D^{- \delta}$, $n \geq 1$, satisfies the strong Doeblin 
condition uniformly in $\varepsilon$  (a similar statement can be found in \cite{FR}, Lemma 3.7.1). 
Therefore, by Lemma~\ref{sdc}, for each $\eta > 0$ and all sufficiently large $n$ (depending on $\eta$), we have the following estimate  for  first term on the right hand side of (\ref{miqu}):
\[
\sup_{x \in \partial D^{-\delta}} | \mathrm{E}  \tau^{ {X}^{x, \varepsilon}_{ \widetilde{\sigma}^{x, \varepsilon}_n},  \varepsilon}(\partial D^{+\delta})  -
\int_{\partial D^{-\delta}} \mathrm{E} \tau^{y, \varepsilon}(\partial D^{+\delta})  \widetilde{\mu}^\varepsilon(dy)  | \leq \eta I(\varepsilon). 
\]
Let us estimate the second term on the right hand side of  (\ref{miqu}). 
By the strong Markov property of the process ${X}^{x, \varepsilon}_t$,
\[
\sup_{x \in \partial D^{-\delta}} | \mathrm{E}  \tau^{x, \varepsilon}(\partial D^{+\delta})  - \mathrm{E}  
 \tau^{{X}^{x, \varepsilon}_{ \widetilde{\sigma}^{x, \varepsilon}_n}, \varepsilon}(\partial D^{+\delta})  |  
\]
\[
\leq \sup_{x \in \partial D^{-\delta}} \mathrm{E} (|   \tau^{x, \varepsilon}(\partial D^{+\delta})  -  
 \tau^{{X}^{x, \varepsilon}_{ \widetilde{\sigma}^{x, \varepsilon}_n}, \varepsilon}(\partial D^{+\delta})  | \chi_{ \{ \widetilde{\sigma}^{x, \varepsilon}_n
\leq \tau^{x, \varepsilon}(\partial D^{+\delta}) \} } ) +
\]
\[
\sup_{x \in \partial D^{-\delta}} \mathrm{E} ( \tau^{x, \varepsilon}(\partial D^{+\delta})   \chi_{ \{ \widetilde{\sigma}^{x, \varepsilon}_n
> \tau^{x, \varepsilon}(\partial D^{+\delta}) \} }  ) + \sup_{x \in \partial D^{-\delta}} \mathrm{E} ( \tau^{{X}^{x, \varepsilon}_{ \widetilde{\sigma}^{x, \varepsilon}_n}, \varepsilon}(\partial D^{+\delta})   \chi_{ \{ \widetilde{\sigma}^{x, \varepsilon}_n
> \tau^{x, \varepsilon}(\partial D^{+\delta}) \} }  )
\]
\[
\leq 2 \sup_{x \in \partial D^{-\delta}}  \mathrm{E} \widetilde{\sigma}^{x, \varepsilon}_n + 
\sup_{x \in \partial D^{-\delta}} \mathrm{P}(\widetilde{\sigma}^{x, \varepsilon}_n > {\tau}^{x, \varepsilon}(\partial D^{+\delta})) 
\sup_{x \in \partial D^{-\delta}} \mathrm{E} {\tau}^{x, \varepsilon}(\partial D^{+\delta}))
\]
\[
\leq 2 \sup_{x \in \partial D^{-\delta}}  \mathrm{E} \widetilde{\sigma}^{x, \varepsilon}_n + 
\sup_{x \in \partial D^{-\delta}} \mathrm{P}(\widetilde{\sigma}^{x, \varepsilon}_n > {\tau}^{x, \varepsilon}(\partial D^{+\delta})) 
 (I(\varepsilon) + J(\varepsilon)).
\]
 From (\ref{ff1}) and Lemma~\ref{jkjl} it follows that for fixed $n$ and all sufficiently small $\varepsilon$ we have
\[
2 \sup_{x \in \partial D^{-\delta}}  \mathrm{E} \widetilde{\sigma}^{x, \varepsilon}_n \leq \eta, ~~~~
\sup_{x \in \partial D^{-\delta}} \mathrm{P}(\widetilde{\sigma}^{x, \varepsilon}_n > {\tau}^{x, \varepsilon}(\partial D^{+\delta}))  \leq \eta.
\]
Therefore, taking absolute value and supremum over $\partial D^{-\delta}$ on both sides of (\ref{miqu}), we obtain
\[
I(\varepsilon) \leq \eta I (\varepsilon) + \eta + \eta (I(\varepsilon) + J(\varepsilon)). 
\]
Since $\eta > 0$ was arbitrary and $J(\varepsilon)$ is bounded from below by a positive constant, this implies that $\lim_{\varepsilon \downarrow 0} I(\varepsilon)/J(\varepsilon) = 0$. From here it follows that
[
\[
\lim_{\varepsilon \downarrow 0} \frac{\mathrm{E} \tau^{x, \varepsilon}(\partial D^{+\delta})}{\int_{\partial D^{-\delta}} \mathrm{E} \tau^{y, \varepsilon}(\partial D^{+\delta}) \widetilde{\mu}^\varepsilon(dy)} = 1
\]
uniformly in $x \in \partial D^{-\delta}$. This, together with (\ref{ff1}), implies (\ref{cdc}).

 From (\ref{cdc}) it follows that 
\begin{equation} \label{frft}
\lim_{\varepsilon \downarrow 0} \frac{\mathrm{E} \tau^{x_1, \varepsilon}(\partial D^{+\delta})}{\mathrm{E} \tau^{x_2, \varepsilon}(\partial D^{+\delta})} = 1
\end{equation}
uniformly in $x_1, x_2 \in \partial D$.
Define two more sequences of stopping times:
\[
\overline{\sigma}^{x, \varepsilon}_0 = 0,
\]
\[
\overline{\tau}^{x, \varepsilon}_n = \inf(t \geq \overline{\sigma}^{x, \varepsilon}_{n-1}: {X}^{x, \varepsilon}_t \in \partial D),~~ n \geq 1,
\]
\[  
\overline{\sigma}^{x, \varepsilon}_n = \inf(t \geq \overline{\tau}^{x, \varepsilon}_n: {X}^{x, \varepsilon}_t \in \partial D^{+ \delta}),~~n \geq 1. 
\]
The discrete-time Markov chains ${X}^{x, \varepsilon}_{ \overline{\tau}^{x, \varepsilon}_n} \in \partial D$ and 
$ {X}^{x, \varepsilon}_{ \overline{\sigma}^{x, \varepsilon}_n} \in \partial D^{+ \delta}$, $n \geq 1$, satisfy the strong 
Doeblin condition for each $\varepsilon > 0$. 
Let $\mu^\varepsilon$ and $\eta^\varepsilon$ be their invariant measures (on $\partial D$ and $\partial D^{+ \delta}$, respectively). 
Since the Lebesgue measure is invariant for the process  ${X}^{x, \varepsilon}_t$,
\[
\frac{{\rm Vol} (D^{+\delta} \setminus D) }{ {\rm Vol} (D) } =  \frac{\int_{\partial D^{+\delta}}  \mathrm{E}  \int_0^{ 
{\tau}^{x, \varepsilon}(\partial D)} 
\chi_{D^{+\delta} \setminus D } ({X}^{x, \varepsilon}_t) dt 
\eta^\varepsilon(dx)}{\int_{\partial D}  \mathrm{E}  \int_0^{ {\tau}^{x, \varepsilon}(\partial D^{+\delta})} 
\chi_{D} ({X}^{x, \varepsilon}_t) dt 
\mu^\varepsilon(dx)},
\]
\[
\frac{{\rm Vol} ( D^{+\delta} \setminus D) }{ {\rm Vol} (D^{+\delta}) } =  \frac{\int_{\partial D^{+\delta}}  \mathrm{E}  \int_0^{ 
{\tau}^{x, \varepsilon}(\partial D)} 
\chi_{D^{+\delta} \setminus D  } ({X}^{x, \varepsilon}_t) dt 
\eta^\varepsilon(dx)}{ \int_{\partial D^{+\delta}}  \mathrm{E}  \int_0^{{\tau}^{x, \varepsilon}(\partial D)} \chi_{D^{+ \delta}} ( {X}^{x, \varepsilon}_t) dt 
 \eta^\varepsilon(dx) + 
 \int_{\partial D} \mathrm{E}  {\tau}^{x, \varepsilon}(\partial D^{+\delta})
\mu^\varepsilon(dx)}.
\]
Therefore, 
\[
\frac{{\rm Vol} (D) }{{\rm Vol} (D^{+\delta} \setminus D) } \int_{\partial D^{+\delta}}  \mathrm{E}  \int_0^{ 
{\tau}^{x, \varepsilon}(\partial D)} 
\chi_{D^{+\delta} \setminus D } ({X}^{x, \varepsilon}_t) dt 
\eta^\varepsilon(dx) \leq
\]
\[
\int_{\partial D} \mathrm{E}  {\tau}^{x, \varepsilon}(\partial D^{+\delta})
\mu^\varepsilon(dx) \leq
\]
\[
\frac{{\rm Vol} (D^{+ \delta}) }{{\rm Vol} (D^{+\delta} \setminus D) } \int_{\partial D^{+\delta}}  \mathrm{E}  \int_0^{
{\tau}^{x, \varepsilon}(\partial D)} 
\chi_{D^{+\delta} \setminus D } ({X}^{x, \varepsilon}_t) dt 
\eta^\varepsilon(dx).
\]
Since the process $ {X}^{x, \varepsilon}_t $  coincides with the Brownian outside $D$, we have
\[
\int_{\partial D^{+\delta}}  \mathrm{E}  \int_0^{
{\tau}^{x, \varepsilon}(\partial D)} 
\chi_{D^{+\delta} \setminus D} ({X}^{x, \varepsilon}_t) dt 
\eta^\varepsilon(dx) = \delta^2(1 + {\alpha}^\varepsilon(\delta)).
\]
Since 
\[
\lim_{\delta \downarrow 0} {\rm Vol} (D^{+ \delta}) = {\rm Vol} (D)
\]
and
\[
\lim_{\delta \downarrow 0} ({\rm Vol} (D^{+\delta} \setminus D)/\delta) = \nu(\partial D),
\]
we obtain that
\begin{equation} \label{ou1}
\int_{\partial D} \mathrm{E}  {\tau}^{x, \varepsilon}(\partial D^{+\delta})
\mu^\varepsilon(dx) = \frac{  {\rm Vol} (D) }{  \nu(\partial D)}(1 + {\alpha}^\varepsilon(\delta)) \delta.
\end{equation}
From (\ref{frft}) and (\ref{ou1}) we obtain that
\[
\sup_{x \in \partial {D} } |
\mathrm{E} \tau^{x, \varepsilon}(\partial D^{+\delta}) - \delta \frac{  {\rm Vol} (D) }{  \nu(\partial D)} | =   \delta {\alpha}^{ \varepsilon}(\delta).
\] 
The statement of the lemma now follows from Lemma~\ref{jkjl} and the strong Markov property of the process.
\qed
\\
\\
{\it Proof of Lemma~\ref{le3rr}.} Recall the definition of the discrete-time Markov chains ${X}^{x, \varepsilon}_{ \overline{\tau}^{x, \varepsilon}_n} \in \partial D$ and 
${X}^{x, \varepsilon}_{ \overline{\sigma}^{x, \varepsilon}_n} \in \partial D^{+ \delta}$, $n \geq 1$,   and their invariant measures $\mu^\varepsilon$ and $\eta^\varepsilon$, introduced
in the proof of the previous lemma. Select $\delta$ sufficiently small so that $x \rightarrow (h, \theta)$ is a diffeomorphism
between $\overline{D^{+ 2 \delta}} \setminus D$ and $[0,2\delta] \times \partial D$. Define
\[
\hat{\varphi}^\delta(x) = \begin{cases} \varphi(\theta(x)), & \mbox{if  } x \in  \overline{D^{+  2 \delta}} \setminus D^{+\delta}  \\ 0, & \mbox{otherwise}. \end{cases}
\]
Since the Lebesgue measure  is invariant for the process ${X}^{x, \varepsilon}_t$, 
\begin{equation} \label{ii0}
\frac{\int_{ D^{+2 \delta}} \hat{\varphi}^\delta(x) dx}{{\rm Vol}(D^{+\delta})} = \frac{ \int_{\partial D} \mathrm{E} \int_0^{\overline{\tau}^{x, \varepsilon}_2} 
  \hat{\varphi}^\delta({X}^{x, \varepsilon}_t) dt \mu^\varepsilon(dx)}{ \int_{\partial D} \mathrm{E} \int_0^{\overline{\tau}^{x, \varepsilon}_2} 
  \chi_{D^{+\delta}} ( {X}^{x, \varepsilon}_t)  dt  \mu^\varepsilon(dx)},
\end{equation}
where, as we recall, $\overline{\tau}^{x, \varepsilon}_2$  is the first time when the process ${X}^{x, \varepsilon}_t$ starting at $x \in \partial D$ returns to $\partial D$ after visiting 
$\partial D^{+\delta}$.  
Observe that
\begin{equation} \label{ii1}
\frac{\int_{D^{+2 \delta}} \hat{\varphi}^\delta(x) dx}{  {\rm Vol}(D^{+\delta})}  = \frac{ \int_{\partial D} \varphi (x)  {\nu (dx)}}{{\rm Vol}(D) }(\delta + o(\delta))~~~{\rm as}~\delta \downarrow 0. 
\end{equation}
From (\ref{ou1}), the strong Markov property of the process, and the fact that there is $C > 0$ such that $\sup_{x \in \partial D^{+\delta}} 
\mathrm{E} \int_0^{{\tau}^{x, \varepsilon}(\partial D)} \chi_{D^{+\delta}} ( {X}^{x, \varepsilon}_t)  dt \leq C \delta^2$,     it follows that
\begin{equation} \label{ii2}
\int_{\partial D} \mathrm{E} \int_0^{\overline{\tau}^{x, \varepsilon}_2} 
  \chi_{D^{+\delta}} ( {X}^{x, \varepsilon}_t)  dt  \mu^\varepsilon(dx) = \frac{  {\rm Vol} (D) }{  \nu(\partial D)}  (1 +  {\alpha}^{\varepsilon}(\delta)  ) \delta.
\end{equation}
By the strong Markov property of the process,
\[
\int_{\partial D} \mathrm{E} \int_0^{\overline{\tau}^{x, \varepsilon}_2} 
  \hat{\varphi}^\delta({X}^{x, \varepsilon}_t) dt \mu^\varepsilon(dx) =  \int_{\partial D}  \int_{\partial D^{+\delta}} \mathrm{E} 
	\int_0^{{\tau}^{y, \varepsilon}(\partial D)} 
  \hat{\varphi}^\delta({X}^{y , \varepsilon}_t) dt \mu^{x, \varepsilon}_\delta (dy) \mu^\varepsilon(dx). 
\]
Using this, and employing (\ref{ii0}), (\ref{ii1}), and (\ref{ii2}) consecutively, we obtain that
\[
\int_{\partial D}  \int_{\partial D^{+\delta}} \mathrm{E} \int_0^{{\tau}^{y, \varepsilon}(\partial D)} 
  \hat{\varphi}^\delta({X}^{y, \varepsilon}_t) dt \mu^{x, \varepsilon}_\delta(dy) \mu^\varepsilon(dx) = 
	\] 
	\[
	\left( \frac{\int_{ D^{+2 \delta}} \hat{\varphi}^\delta(x) dx}{{\rm Vol}(D^{+\delta})} \right) 
	{ \int_{\partial D} \mathrm{E} \int_0^{\overline{\tau}^{x, \varepsilon}_2} 
  \chi_{D^{+\delta}} ( {X}^{x, \varepsilon}_t)  dt  \mu^\varepsilon(dx)} =
	\]
	\[
\frac{ \int_{\partial D} \varphi (x)  {\nu (dx)}}{{\rm Vol}(D) }(\delta + o(\delta))	{ \int_{\partial D} \mathrm{E} \int_0^{\overline{\tau}^{x, \varepsilon}_2} 
  \chi_{D^{+\delta}} ( {X}^{x, \varepsilon}_t)  dt  \mu^\varepsilon(dx)} =
	\]
\[
\frac{\delta^2(1 + \alpha^\varepsilon(\delta))}{  \nu(\partial D)} { \int_{\partial D} \varphi (x)  {\nu}(dx) } = 
	\delta^2(1 + \alpha^\varepsilon(\delta)) { \int_{\partial D} \varphi (x) \overline{\nu} (dx)}.
\]
Observe that 
\[
\sup_{y \in \partial D^{+\delta}} |
\mathrm{E} \int_0^{{\tau}^{y, \varepsilon}(\partial D)} 
  \hat{\varphi}^\delta({X}^{y, \varepsilon}_t) dt  - \varphi(\theta(y)) \delta^2 |  = \delta^2 \alpha^{\varepsilon}(\delta). 
\]
Therefore,
\[
\int_{\partial D}  \int_{\partial D^{+\delta}} \varphi(\theta(y))  \mu^{x, \varepsilon}_\delta(dy) \mu^\varepsilon(dx) = 
 (1 + \alpha^\varepsilon(\delta)) { \int_{\partial D} \varphi (x)  \overline{\nu}(dx)}.
\]
Finally, we observe that $\mu^{x, \varepsilon}_\delta(y)$ asymptotically (as $\varepsilon \downarrow 0$) does not depend on $x$, i.e., 
\[
\lim_{\varepsilon \downarrow 0} \frac{ \int_{\partial D^{+\delta}} \varphi(\theta(y)) \mu^{x, \varepsilon}_\delta(dy)}{\int_{\partial D}  \int_{\partial D^{+\delta}} \varphi(\theta(y)) \mu^{z, \varepsilon}_\delta(dy) \mu^\varepsilon(dz)} = 1
\]
uniformly in $x \in \partial D$. 
The proof of this statement is the same as that of (\ref{cdc}). We conclude that
\[
\sup_{x \in \partial D} |
\int_{\partial D^{+\delta}} \varphi(\theta(y)) \mu^{x, \varepsilon}_\delta(dy) - { \int_{\partial D} \varphi(y) \overline{\nu} (dy)}| =  \alpha^{\varepsilon}(\delta),
\]
which completes the proof of the lemma.
\qed
\section{A non-standard boundary problem}
In this section, we show that the expected time for the limiting process to reach a given set is the solution to a non-standard boundary problem
for an elliptic PDE. Similar elliptic and parabolic problems come up when considering other quantities associated with the limiting process (probability of reacihing
one given set prior to another, expected time spent by the process in a given set prior to time $t$, etc.). We restrict ourselves to one example.  

Let $D_1,...,D_n \subset \mathbb{T}^d $ be as in Section~\ref{intro}, and let $F $ be
a domain with $C^3$-smooth boundary such that $\overline{F} \subset   U = \mathbb{T}^d \setminus \bigcup_{k=1}^n \overline{D}_k$. For $x \in \mathbb{T}^d \setminus F$,
 let
\[
\check{\tau}^{x}(\overline{F}) =  \inf \{t\geq 0: Y^{\varphi(x)}_t \in \overline{F} \},~~~~u(x) = \mathrm{E} \check{\tau}^{x}(\overline{F}),
\]
where $Y^{\varphi(x)}_t$ is the process constructed in Section~\ref{dttp}. It is clear that $\mathrm{E} \check{\tau}^{x}(\overline{F})$ is constant on each of the
sets $D_1,...,D_n$.  Our goal is to express $\mathrm{E} \check{\tau}^{x}(\overline{F})$
as a solution of an elliptic problem on $\overline{U} \setminus F$. 

Let $u \in C^2(\overline{U} \setminus F)$ solve the following problem:
\[
\frac{\Delta {u}(x)}{2} = -1,~~~~x \in U;
\]
\[
u(x) = c_k,~~x \in \partial D_k,~~k=1,...,n;
\]
\[
u(x) = 0,~~x \in \partial F;
\]
\[
\int_{\partial D_k} \langle \nabla u(x) , n(x) \rangle \nu_k(d x) - {\rm Vol}(D_k) = 0,~~~k =1,...,n,
\]
where  $n(x)$ is the unit exterior  normal at $x \in \partial D_k$ (external with respect to $U$) and $\nu_1,...,\nu_n$ are the  $(d-1)$-dimensional  volume measures on $\partial D_1,...,\partial D_n$, respectively. The constants $c_k$, $k =1,...,n$, are not prescribed, i.e., solving the problem includes finding the constants. 
The solution $u \in C^2(\overline{U} \setminus F)$ exists and is unique, as can be easily justified with the same arguments as those used in Section~\ref{dttp}
to show the existence and uniqueness of a similar elliptic problem.

Let $\tilde{u} \in C(\mathbb{T}^d) \bigcap C^2(\overline{U})$ coinside with $u$ on $\overline{U} \setminus F$ and be constant on each of the sets $D_1,...,D_n$ (i.e.,
we extend $u$ as constants in  $D_1,...,D_n$ and as a $C^2$ function in $F$).
\begin{lemma} For $x \in \mathbb{T}^d \setminus F$, we have 
$\mathrm{E} \check{\tau}^{x}(\overline{F}) = \tilde{u}(x)$. 
\end{lemma}
\noindent
{\it Sketch of the proof.} From the convergence of the process $Y^{x, \varepsilon}_t = \varphi(X^{x,\varepsilon}_t)$ 
to the process $Y^{\varphi(x)}_t$ (Theorem~\ref{mt1a}), it easily follows that
\[
\mathrm{E} \check{\tau}^{x}(\overline{F}) = \lim_{\varepsilon \downarrow 0} \mathrm{E} {\tau}^{x, \varepsilon}(\overline{F}),
\]
where $\tau^{x, \varepsilon}(\overline{F})$ is the first time when the process  $X^{x,\varepsilon}_t$ reaches $\overline{F}$. We also
note that Proposition~\ref{mainpro} remains valid with the stopping time ${\tau}^{x, \varepsilon}(\overline{F})$ instead of a fixed time $t$
(the proof of the proposition requires only minor modifications). Applying this to $f \in \mathcal{D}(A)$ defined 
as $f(x) =  \tilde{u}(\varphi^{-1}(x))$, $x \in U'$, we obtain
\[
\lim_{\varepsilon \downarrow 0} \left(-\tilde{u}(x) + \mathrm{E} {\tau}^{x, \varepsilon}(\overline{F}) \right) = \lim_{\varepsilon \downarrow 0}
 \mathrm{E} \left(\tilde{u}(X^{x, \varepsilon}_{{\tau}^{x, \varepsilon}(\overline{F}) }) -\tilde{u}(x) + \mathrm{E} {\tau}^{x, \varepsilon}(\overline{F}) \right) = 
\]
\[
\lim_{\varepsilon \downarrow 0} \mathrm{E} \left(f(Y^{x, \varepsilon}_{{\tau}^{x, \varepsilon}(\overline{F}) }) - f(x) - 
 \frac{1}{2} \int_0^{{\tau}^{x, \varepsilon}(\overline{F}) } \Delta f (Y^{x, \varepsilon}_u) du \right) = 0,~~~~x \in \mathbb{T}^d \setminus F.
\]
(As before, we used the convention that a function $f$ defined on $U'$ can also be viewed as a function on $\mathbb{T}^d$.) Thus
\[
\mathrm{E} \check{\tau}^{x}(\overline{F}) = \lim_{\varepsilon \downarrow 0} \mathrm{E} {\tau}^{x, \varepsilon}(\overline{F}) = \tilde{u}(x).
\]
\qed
\\
\\
\noindent {\bf \large Acknowledgements}: While working on this
article, M. Freidlin was supported by NSF grant DMS-1411866
and L. Koralov was supported by NSF grant DMS-1309084 and ARO grant W911NF1710419.
\\
\\

\end{document}